\documentclass[11pt,reqno]{amsart}
\usepackage{amsmath,amssymb,graphicx}
\newtheorem{thm}{Theorem}[section]
\newtheorem{cor}[thm]{Corollary}
\newtheorem{lem}[thm]{Lemma}
\newtheorem{prop}[thm]{Proposition}
\theoremstyle{definition}

\theoremstyle{remark}
\newtheorem{rem}[thm]{Remark}
\theoremstyle{remark}

\numberwithin{equation}{section}

\begin{document}

\title[A study of $BES_{0}(\delta)$ processes
$(0<\delta<2)$]{Some random times and martingales associated with
$BES_{0}(\delta)$ processes $(0<\delta<2)$}
\author{Ashkan Nikeghbali}
\address{Laboratoire de Probabilit\'es et Mod\`eles Al\'eatoires, Universit\'e Pierre et Marie Curie, et CNRS UMR 7599,
175 rue du Chevaleret F-75013 Paris, France.} \curraddr{American
Institute of Mathematics 360 Portage Ave Palo Alto, CA 94306-2244
USA} \email{ashkan@aimath.org}

\subjclass[2000]{05C38, 15A15; 05A15, 15A18} \keywords{Bessel
processes, Random Times, General theory of stochastic processes,
Optional stopping theorem, Martingales, zeros of Bessel processes.}

\date{\today}
\begin{abstract}
In this paper, we study Bessel processes of dimension
$\delta\equiv2(1-\mu)$, with $0<\delta<2$, and some related
martingales and random times. Our approach is based on martingale
techniques and the general theory of stochastic processes (unlike
the usual approach based on excursion theory), although for
$0<\delta<1$, these processes are even not semimartingales. The last
time before $1$ when a Bessel process hits $0$, called $g_{\mu}$,
plays a key role in our study: we characterize its conditional
distribution and extend Paul L\'{e}vy's arc sine law and a related
result of Jeulin about the standard Brownian Motion. We also
introduce some remarkable families of martingales related to the
Bessel process, thus obtaining in some cases a one parameter
extension of some results of Az\'{e}ma and Yor in the Brownian
setting: martingales which have the same set of zeros as the Bessel
process and which satisfy the stopping theorem for $g_{\mu}$, a one
parameter extension of Az\'{e}ma's second martingale, etc.
Throughout our study, the local time of the Bessel process also
plays a central role and we shall establish some of its elementary
properties.
\end{abstract}
\maketitle
\section{Introduction}
Bessel processes afford basic examples of diffusion processes and
their systematic study  was initiated in McKean \cite{mckean}. There
exist many results for Bessel processes of dimension $\delta>2$,
obtained thanks to stochastic calculus and martingale techniques;
when $\delta\in(0,2)$, the same methods do not apply anymore:
indeed, for $\delta\in(0,1)$, Bessel processes are even not
semimartingales (see \cite{revuzyor}, chapter XI for more details
and references).

The main tool to study Bessel processes of dimension
$\delta\in(0,2)$ is the powerful theory of excursions of Markov
processes, as developed e.g. by Getoor (\cite{Getoor}) and
previously initiated by It\^{o} (\cite{ito}). Among other works
related to the study of Bessel processes, we can cite the seminal
work of Bertoin (\cite{bertoin}) who developed an excursion theory
for the Bessel process of dimension $\delta\in(0,1)$ and its drift
term and the work of Barlow, Pitman and Yor (\cite{barlowpitmanyor})
who proved a multidimensional extension of the arc sine law, for
Bessel processes of dimension $\delta\in(0,2)$, using excursion
theory. We should also mention here the work of Az\'{e}ma
\cite{azemaferme}, who developed a theory of excursions for closed
optional sets, which led him to discover the so called first and
second Az\'{e}ma's martingales (the reader can refer to
\cite{delmaismey} for an introduction and more references on this
topic). This approach is interesting when, for example, one studies
martingales or submartingales in the filtration of the zeros of a
diffusion process.

In this paper, we shall consider the couple
$\left(R_{t},L_{t}\right)$, where $\left(R_{t}\right)$ is a Bessel
process of dimension $\delta\equiv2(1-\mu)\in(0,2)$, starting from
$0$, and $\left(L_{t}\right)$ a choice of its local time at level
zero. We shall associate with $R$ the honest time:
$$g_{\mu}\left(T\right)\equiv\sup\left\{u\leq
T:\;R_{u}=0\right\},$$where $T>0$ is a fixed time (we shall simply
note $g_{\mu}$ when $T=1$). The case $\mu=\dfrac{1}{2}$ has received
much attention in the literature: the distribution of $g_{1/2}$ was
obtained by L\'{e}vy (\cite{levy}), the supermartingale
$Z_{t}\equiv\mathbb{P}\left[g_{1/2}>t\mid\mathcal{F}_{t}\right]$ and
the dual predictable projection of $\mathbf{1}_{\left(g_{1/2}\leq
t\right)}$, which play a key role in the general theory of
stochastic processes (see \cite{azema}), were obtained by Jeulin
(\cite{jeulin}). Moreover, the latter quantities computed by Jeulin
have revealed to be useful on the one hand to Marc Yor in the study
of martingales which have the same zeros as the standard Brownian
Motion (see \cite{zurich}, \cite{columbia}), and on the other hand
to Az\'{e}ma et alii. (\cite{azemajeulinknightyor}) in the study of
the stopping theorem when stopping times are replaced with $g_{1/2}$
(see also \cite{Ashkannonstop}). The computations performed by
Jeulin have also been very useful in the mathematical models of
default times (see \cite{elliotjeanbyor} for examples and more
references). This paper has three main aims:
\begin{itemize}
\item to extend the computations of Jeulin to the case of Bessel
processes of dimension $\delta\in(0,2)$ and then to generalize the
above mentioned results of Yor;
\item to illustrate the results in \cite{azemajeulinknightyor} and
\cite{Ashkannonstop} on the stopping theorem when considered at an
honest time and more generally to give a larger class of examples
than the usual example of the standard Brownian Motion;
\item to illustrate martingale techniques in a setting where only
excursion theory is used (except for a forthcoming paper by
Roynette, Vallois and Yor \cite{roynete} which deals with the
penalization of Bessel process of dimension $0<\delta<2$ with a
function of its local time at $0$).
\end{itemize}
More precisely, the paper is organized as follows:

\noindent In Section 2, we recall and establish some basic facts
about the local time of a Bessel process and then we compute
explicitly the supermartingale
$Z_{t}^{g_{\mu}}\equiv\mathbb{P}\left[g_{\mu}>t\mid\mathcal{F}_{t}\right]$
and the dual predictable projection of
$\mathbf{1}_{\left(g_{\mu}\leq t\right)}$. We then use them to
obtain a one parameter extension of L\'{e}vy's arc sine law. We thus
recover a result which Barlow, Pitman and Yor
(\cite{barlowpitmanyor}) have found using excursion theory (in fact,
this result is due to Dynkin \cite{dynkin} by completely different
means).

\noindent In Section 3, we investigate further for the conditional
distribution of $g_{\mu}$ and then test the stopping theorem on
martingales of the form:
$$M_{t}^{h}\equiv\mathbb{E}\left[h\left(g_{\mu}\right)\mid\mathcal{F}_{t}\right].$$More
precisely, we show how
$\mathbb{E}\left[M_{\infty}^{h}\mid\mathcal{F}_{g_{\mu}}\right]$ and
$M_{g_{\mu}}^{h}$ differ (in particular we recover the results in
\cite{zurich} when we take $\mu=\dfrac{1}{2}$).

\noindent In Section 4, we characterize the martingales
$\left(M_{t}\right)$ which have the same set of zeros as $R$ and use
this characterization to give many examples of martingales which
satisfy the stopping theorem with $g_{\mu}$, i.e.:
$$\mathbb{E}\left[M_{\infty}\mid\mathcal{F}_{g_{\mu}}\right]=M_{g_{\mu}}^{h}.$$
Here again, we obtain some natural extensions of the results of Yor
in the Brownian setting (\cite{zurich}, \cite{columbia}).

\noindent Eventually, in Section 5, using some not so well known
results of Yor (\cite{zurichI}) about Bessel meanders of dimension
$0<\delta<2$, we obtain a one dimensional extension of Az\'{e}ma's
second martingale. We compare our result with some results of
Az\'{e}ma (\cite{azemaferme}) and C. Rainer (\cite{rainer}) about
projections of a diffusion on its slow filtration to recover the
L\'{e}vy measure of the zeros of $R$. More generally, we give many
examples of $\left(\mathcal{F}_{g_{\mu}\left(t\right)}\right)$
martingales by computing the projections of some carefully selected
martingales, reminiscent of the Az\'{e}ma-Yor martingales . We shall
also combine these arguments with Doob's maximal identity to obtain
some local time estimates; more precisely, inspired by the work of
Knight \cite{Knight73,knight78}, and some recent lecture given by
Marc Yor at Columbia University (\cite{columbia}), we compute
explicitly the following probabilities:
$$\mathbb{P}\left(\exists t\geq
0,\;R_{t}>\varphi\left(L_{t}\right)\right),$$and
$$\mathbb{P}\left(\exists t\leq
\tau_{u},\;R_{t}>\varphi\left(L_{t}\right)\right),$$where $\varphi$
is a positive Borel function and $\left(\tau_{u}\right)$ the
right-continuous inverse of the local time $\left(L_{u}\right)$.
\section{The local time of the Bessel process and an extension of the arc sine law}
\subsection{Basic facts about the local time of a Bessel
process}We shall now precisely define what we mean by the local time
of a Bessel process: indeed, in the literature, one can find
different normalizations for the local time (see the forthcoming
work \cite{donatimartin}). Our approach is based on a result of
Biane and Yor about powers of Bessel processes (\cite{bianeyor}).

Let $\left( \Omega ,\mathcal{F},\left( \mathcal{F}_{t}\right),%
\mathbb{P}\right) $ be a filtered probability space, where the filtration $%
\left( \mathcal{F}_{t}\right) $\ is generated by a Bessel process
$\left( R_{t}\right) _{t\geq 0}$, of dimension $\delta\in(0,2)$,
starting from $0$. We associate with $\delta$ two other parameters
$\nu$ and $\mu$:
$$\nu\equiv\dfrac{\delta}{2}-1,\;\mu\equiv-\nu.$$ We note that
$\nu\in(-1,0)$ and $\mu\in(0,1)$. We first recall some basic facts
about Bessel processes of index $\nu\in(-1,0)$ (see for example
\cite{borodinsalminen}, \cite{itomackean}, \cite{revuzyor}). The
process $\left( R_{t}\right)$ is an $\mathbb{R}_{+}$-valued
diffusion whose infinitesimal generator is defined by:
$$\mathcal{L}f\left(r\right)=\dfrac{1}{2}\dfrac{d^{2}f}{dr^{2}}+\dfrac{1-2\mu}{2r}\dfrac{df}{dr},$$on
the domain:
$$\mathcal{D}=\left\{f:\mathbb{R}_{+}\rightarrow\mathbb{R};\;\mathcal{L}f\in\mathcal{C}_{b}\left(\mathbb{R}_{+}\right),\;\lim_{r\rightarrow0}r^{1-2\mu}f'\left(r\right)=0\right\}.$$
$\left(R_{t}\right)$ is a recurrent diffusion, and
$\left\{0\right\}$ is a reflecting point. Its scale function is
given by:
$$s\left(x\right)=x^{-2\nu}(=x^{2\mu}),$$ and its speed measure by:
$$m_{\nu}\left(dx\right)=\dfrac{x^{2\nu+1}}{\left|\nu\right|}dx(=\dfrac{x^{1-2\mu}}{\mu}dx).$$ The semigroup
of $\left(R_{t}\right)$, with respect to the speed measure, is given
by:
\begin{eqnarray*}
  p^{\left(\nu\right)}\left(t;x,y\right) &=& \dfrac{\left|\nu\right|}{t}\exp\left(-\dfrac{x^{2}+y^{2}}{2t}\right)I_{\nu}\left(xy\right)\left(xy\right)^{-\nu},\;x>0 \\
  p^{\left(\nu\right)}\left(t;0,y\right) &=&
  \dfrac{\left|\nu\right|}{2^{\nu}t^{\nu+1}\Gamma\left(\nu+1\right)}\exp\left(-\dfrac{y^{2}}{2t}\right),
\end{eqnarray*}where $I_{\nu}\left(x\right)=\sum_{k=0}^{\infty}
\frac{\left(x/2\right)^{\nu+2k}}{k!\Gamma\left(\nu+k+1\right)}$ is
the modified Bessel function  and $\Gamma$ is Euler's gamma
function.

Now, we shall look at a conveniently chosen power of the Bessel
process to define the local time at $0$ of $R$:
\begin{prop}[and Definition]\label{bianyorpuissance}
There exists a reflected Brownian Motion $\left(\gamma_{t}\right)$,
on the same probability space, such that:
\begin{equation}\label{puissancedebessel}
    R_{t}^{2\mu}=2\mu\gamma_{\int_{0}^{t}R_{u}^{2\left(2\mu-1\right)}du}.
\end{equation}Inspired by the Brownian case ($\mu=\frac{1}{2}$), we
shall take as a definition for $\left(L_{t}\right)$, the local time
at $0$ of $\left( R_{t}\right)$, the unique increasing process
$\left(L_{t}\right)$ such that $$N_{t}\equiv R_{t}^{2\mu}-L_{t},$$
is a martingale. Moreover, we have: $\langle
N\rangle_{t}=4\mu^{2}\int_{0}^{t}duR^{2\left(2\mu-1\right)}_{u}$,
and
$L_{t}=\ell_{4\mu^{2}\int_{0}^{t}R_{u}^{2\left(2\mu-1\right)}du}$,
where $\left(\ell_{u}\right)$ is  chosen such that
$\left(\gamma_{u}-\ell_{u}\right)_{u\geq 0}$ is an
$\left(\sigma\left\{\gamma_{s},\;s\leq u\right\}\right)_{u\geq 0}$
martingale.
\end{prop}
\begin{proof}
It suffices to prove (\ref{puissancedebessel}); it is a consequence
of a result of Biane and Yor about powers of Bessel processes
(\cite{bianeyor}, or Proposition 1.11 p.447 in \cite{revuzyor}): if
$R$ is a Bessel process of index $\nu$, and if $p$ and $q$ are two
conjugate numbers such that $\frac{1}{p}+\frac{1}{q}=1$, with
$\nu>-\frac{1}{q}$, then there exists a Bessel process
$\widehat{R}$, of index $q\nu$, such that:
$$qR_{u}^{1/q}=\widehat{R}_{\int_{0}^{t}duR_{u}^{-2/p}}.$$
\end{proof}
One should be very careful with the choice or the normalization of
the local time: indeed, from the general theory of local times for
diffusion processes (see for example \cite{borodinsalminen}), there
exists a jointly continuous family $\left(L_{t}^{x};\;x\geq0,
t\geq0\right)$, such that the following occupation formula holds:
\begin{equation}\label{formuleoccupationgen}
\int_{0}^{t}h\left(R_{u}\right)du=C\int_{0}^{\infty}h\left(x\right)L_{t}^{x}x^{1-2\mu}dx,
\end{equation}
for every Borel function
$h:\mathbb{R}_{+}\rightarrow\mathbb{R}_{+}$. The choice of
$L_{t}\equiv L_{t}^{0}$ determines the constant $C$ (see
\cite{donatimartin} for a detailed discussion about the different
normalizations found in the literature). More precisely:
\begin{prop}\label{costeoccupation}
Let $h:\mathbb{R}_{+}\rightarrow\mathbb{R}_{+}$ be a Borel function.
Then, with our choice for $\left(L_{t}\right)$, the following
occupation formula holds:
\begin{eqnarray}\label{formuleoccupationbessel}
    \int_{0}^{t}h\left(R_{u}\right)du &=&
    \dfrac{1}{\mu}\int_{0}^{\infty}h\left(x\right)L_{t}^{x}x^{1-2\mu}dx;
    \\
     &=& \int_{0}^{\infty}h\left(x\right)L_{t}^{x}m\left(dx\right).
\end{eqnarray}Consequently, $C=\dfrac{1}{\mu}$.
\end{prop}
\begin{proof}
Taking the expectation of both sides in (\ref{formuleoccupationgen})
yields:
$$\dfrac{1}{\mu}\int_{0}^{t}dup^{(\nu)}\left(u;0,x\right)=C\mathbb{E}\left[L_{t}^{x}\right].$$Now,
letting $x\rightarrow0$ in the above equation, we obtain:
$$\dfrac{2^{\mu}t^{\mu}}{\mu\Gamma\left(1-\mu\right)}=C\mathbb{E}\left[R_{t}^{2\mu}\right].$$But
since
$$R_{t}^{2}\sim 2t\;gam(1-\mu),$$where $gam(1-\mu)$ denotes a
standard Gamma variable of parameter $1-\mu$, it follows that:
$$C=\dfrac{1}{\mu}.$$
\end{proof}
\begin{rem}\label{remarquetpslocal}
We could obtain the occupation formula without using general results
about diffusion processes, but just using Proposition
\ref{bianyorpuissance}, the occupation formula for the standard
Brownian Motion and time change arguments. Indeed, for any Borel
function $f:\mathbb{R}_{+}\rightarrow\mathbb{R}_{+}$, we have:
$$\int_{0}^{t}f\left(R_{u}\right)du=\int_{0}^{t}g\left(\gamma_{A_{u}}\right)du,$$where
$\gamma$ denotes a reflected Brownian Motion, $A_{t}\equiv\langle
N\rangle_{t}=4\mu^{2}\int_{0}^{t}duR^{2\left(2\mu-1\right)}_{u}$ and
$g\left(x\right)\equiv f\left(x^{1/2\mu}\right)$. If we note
$\left(\rho_{t}\right)$ the right continuous inverse of
$\left(A_{t}\right)$, we have:
$$\int_{0}^{t}g\left(\gamma_{A_{u}}\right)du=\int_{0}^{t}g\left(\gamma_{A_{u}}\right)\dfrac{R_{u}^{2\left(1-2\mu\right)}}{4\mu^{2}}dA_{u}=
\int_{0}^{A_{t}}g\left(\gamma_{u}\right)\dfrac{R_{\rho_{u}}^{2\left(1-2\mu\right)}}{4\mu^{2}}du.$$
But $R_{\rho_{u}}^{2\mu}=\gamma_{u}$, hence:
$$\int_{0}^{t}g\left(\gamma_{A_{u}}\right)du=\dfrac{1}{4\mu^{2}}\int_{0}^{A_{t}}g\left(\gamma_{u}\right)\gamma_{u}^{\frac{1-2\mu}{\mu}}du.$$
Now, the occupation formula for the reflected Brownian Motion
yields:
$$\dfrac{1}{4\mu^{2}}\int_{0}^{A_{t}}g\left(\gamma_{u}\right)\gamma_{u}^{\frac{1-2\mu}{\mu}}du=
\dfrac{1}{4\mu^{2}}\int_{0}^{\infty}f\left(x^{1/2\mu}\right)x^{\frac{1-2\mu}{\mu}}\lambda_{A_{t}}^{x}dx,$$where
$\left(\lambda_{t}^{x},\;t\geq0,x\geq0\right)$ denotes the family of
local times of $\gamma$. A  straightforward change of variables in
the previous integral yields:
$$\dfrac{1}{4\mu^{2}}\int_{0}^{\infty}f\left(x^{1/2\mu}\right)x^{\frac{1-2\mu}{\mu}}\lambda_{A_{t}}^{x}dx=
\dfrac{1}{2\mu}\int_{0}^{\infty}f\left(x\right)x^{1-2\mu}\lambda_{A_{t}}^{x^{2\mu}}dx.$$
Now, plugging all these informations together, we obtain:
$$\int_{0}^{t}f\left(R_{u}\right)du=\dfrac{1}{2\mu}\int_{0}^{\infty}f\left(x\right)x^{1-2\mu}\lambda_{A_{t}}^{x^{2\mu}}dx.$$
Hence,
$$\int_{0}^{t}f\left(R_{u}\right)du=\dfrac{1}{\mu}\int_{0}^{\infty}f\left(x\right)x^{1-2\mu}L_{t}^{x}dx,$$where
$$L_{t}^{x}=\dfrac{\lambda_{A_{t}}^{x^{2\mu}}}{2}.$$This result is
consistent with our choice for the local time since
$\lambda_{t}^{0}=2\ell_{t}$. Consequently, with this time change
method, we have an explicit expression for the local time in both
variables $\left(t,x\right)$, with respect to the local time of the
reflected Brownian Motion. This completes the result in Proposition
\ref{bianyorpuissance}.
\end{rem}
\begin{rem}
The occupation formula and the scaling property for the Bessel
processes (Bessel processes have the Brownian scaling property,
\cite{revuzyor}, chapter XI) entail that $L_{t}$ is distributed as
$t^{\mu}L_{1}$.
\end{rem}
We now give as a corollary of Proposition \ref{costeoccupation} or
Remark \ref{remarquetpslocal} a limit theorem for some integrals of
Bessel processes:
\begin{cor}
If $f$ is a Borel function such that
$$\int_{0}^{\infty}dx\left|f\left(x\right)\right|x^{1-2\mu}<\infty,$$then:
$$\lim_{n\rightarrow\infty}n^{\delta}\int_{0}^{\cdot}f\left(nR_{u}\right)du=\left(\dfrac{1}{\mu}\int_{0}^{\infty}dxf\left(x\right)x^{1-2\mu}\right)L\;a.s.$$
\end{cor}
\begin{proof}
By the occupation formula, for every $t\geq0$:
$$n^{\delta}\int_{0}^{t}f\left(nR_{u}\right)du=\dfrac{1}{\mu}\int_{0}^{\infty}f\left(x\right)L_{t}^{x/n}x^{1-2\mu}dx.$$
Now, for fixed $t$, the result follows from the fact that $x\mapsto
L_{t}^{x}$ is $a.s.$ continuous and has compact support and from the
dominated convergence theorem. The result holds simultaneously for
every rational $t$, and if $f$ is positive, the result follows by
increasing limits. In the general case, it suffices to decompose $f$
as the difference of its positive and negative parts.
\end{proof}
We shall also need the following extension of the occupation times
formula:
\begin{lem}\label{occupationgeneralisee}
Let $h:\mathbb{R}_{+}\times\mathbb{R}_{+}\rightarrow\mathbb{R}_{+}$
be a Borel function. Then, almost-surely,
$$\int_{0}^{t}du\;h\left(u,R_{u}\right)=\dfrac{1}{\mu}\int_{0}^{\infty}da\;a^{1-2\mu}\int_{0}^{t}dL_{u}^{a}h\left(u,a\right).$$
\end{lem}
\begin{proof}
The formula is easily checked for $h(u,a)\equiv
\mathbf{1}_{\left]\alpha,\beta\right]}(u)f(a)$, and the general
result follows by a monotone class argument.
\end{proof}
\subsection{A one parameter extension of the arc sine law}
 Now, define:
$$g_{\mu}\equiv\sup\left\{t\leq1:\;R_{t}=0\right\},$$and more
generally, for $T>0$, a fixed time,
$$g_{\mu}\left(T\right)\equiv\sup\left\{t\leq T:\;R_{t}=0\right\}.$$
\begin{prop}\label{surmartdegamma}
Let $\mu\in(0,1)$, and let $\left(R_{t}\right)$ be a Bessel process
of dimension $\delta=2\left(1-\mu\right)$. Then, we have:
\begin{enumerate}
\item
$$Z_{t}^{g_{\mu}}\equiv\mathbb{P}\left[g_{\mu}>t\mid\mathcal{F}_{t}\right]=\dfrac{1}{2^{\mu-1}\Gamma\left(\mu\right)}\int_{\frac{R_{t}}{\sqrt{1-t}}}^{\infty}
dyy^{2\mu-1}\exp\left(-\frac{y^{2}}{2}\right);$$
\item The dual predictable projection $A_{t}^{g_{\mu}}$of $\mathbf{1}_{\left(g_{\mu}\leq
t\right)}$ is:
$$A_{t}^{g_{\mu}}=\frac{1}{2^{\mu}\Gamma\left(1+\mu\right)}\int_{0}^{t\wedge
1}\dfrac{dL_{u}}{\left(1-u\right)^{\mu}},$$i.e. for every
nonnegative predictable process $\left(x_{t}\right)$,
$$\mathbb{E}\left[x_{g_{\mu}}\right]=\dfrac{1}{2^{\mu}\Gamma\left(1+\mu\right)}\mathbb{E}\left[\int_{0}^{1}dL_{u}\dfrac{x_{u}}{\left(1-u\right)^{\mu}}\right].$$
\end{enumerate}
\end{prop}
\begin{rem}
Some straightforward change of variables allows us to rewrite
$Z_{t}^{g_{\mu}}$ as:
\begin{eqnarray}\label{autreformepourzt}
  Z_{t}^{g_{\mu}} &=& \dfrac{1}{\Gamma\left(\mu\right)}\int_{\frac{R_{t}^{2}}{2\left(1-t\right)}}^{\infty}dz\;z^{\mu-1}\exp\left(-z\right) \\
   &=& \dfrac{R_{t}^{2\mu}}{2^{\mu}\Gamma\left(\mu\right)\left(1-t\right)^{\mu}}\int_{1}^{\infty}dz\;z^{\mu-1}\exp\left(-\dfrac{zR_{t}^{2}}{2\left(1-t\right)}\right).
\end{eqnarray}
\end{rem}
\begin{proof}
$(1).$ We have:
$$Z_{t}^{g_{\mu}}=1-\mathbb{P}\left[g_{\mu}\leq
t\mid\mathcal{F}_{t}\right]=1-\mathbb{P}\left[d_{t}>
1\mid\mathcal{F}_{t}\right],$$where (using the Markov property),
$$d_{t}=\inf\left\{u\geq t;\;R_{u}=0\right\}=t+\inf\left\{u\geq
0;\;R_{t+u}=0\right\}=t+\theta\circ H_{0},$$where $\theta$ is the
shift operator and $H_{0}$ is the first hitting time of $0$.
Consequently, using the Markov property, we obtain:
\begin{equation}\label{aqsf}
Z_{t}^{g_{\mu}}=1-\mathbb{P}_{R_{t}}\left[H_{0}> 1-t\right].
\end{equation}
 Now, following Borodin and Salminen (p.
70-71), if for $-\nu>0$, $\mathbb{P}_{0}^{(-\nu)}$ denotes the law
of a Bessel process of parameter $-\nu$, starting from $0$, then the
law of $L_{y}\equiv\sup\left\{t:\;R_{t}=y\right\}$, is given by:
$$\mathbb{P}_{0}^{(-\nu)}\left(L_{y}\in
dt\right)=\dfrac{y^{-2\nu}}{2^{-\nu}\Gamma\left(-\nu\right)t^{-\nu+1}}\exp\left(-\frac{y^{2}}{2t}\right)dt.$$Now,
from the time reversal property for Bessel processes
(\cite{borodinsalminen} p.70), we have:
$$\mathbb{P}_{x}\left[H_{0}\in dt\right]=\mathbb{P}_{0}^{(-\nu)}\left(L_{x}\in
dt\right);$$consequently, from (\ref{aqsf}), we have (recall
$\mu=-\nu$):
$$Z_{t}^{g_{\mu}}=1-\dfrac{R_{t}^{2\mu}}{2^{\mu}\Gamma\left(\mu\right)}\int_{1-t}^{\infty}du\frac{\exp\left(-\frac{R_{t}^{2}}{2u}\right)}{u^{1+\mu}},$$and
the desired result is obtained by straightforward change of
variables in the above integral.

$(2)$ is a consequence of It\^{o}'s formula applied to
$Z_{t}^{g_{\mu}}$ and the fact that $N_{t}\equiv R_{t}^{2\mu}-L_{t}$
is a martingale and $\left(dL_{t}\right)$ is carried by
$\left\{t:\;R_{t}=0\right\}$.
\end{proof}
\begin{rem}
The previous proof can be applied mutatis mutandis to obtain:
$$\mathbb{P}\left[g_{\mu}(T)>t\mid\mathcal{F}_{t}\right]=\dfrac{1}{2^{\mu-1}\Gamma\left(\mu\right)}\int_{\frac{R_{t}}{\sqrt{T-t}}}^{\infty}
dyy^{2\mu-1}\exp\left(-\frac{y^{2}}{2}\right);$$and
$$A_{t}^{g_{\mu}(T)}=\frac{1}{2^{\mu}\Gamma\left(1+\mu\right)}\int_{0}^{t\wedge
T}\dfrac{dL_{u}}{\left(T-u\right)^{\mu}}.$$
\end{rem}
When $\mu=\dfrac{1}{2}$, $R_{t}$ can be viewed as
$\left|B_{t}\right|$, the absolute value of a standard Brownian
Motion. Thus, we recover as a particular case of our framework the
celebrated example of the last zero before $1$ of a standard
Brownian Motion (see \cite{jeulin} p.124, or \cite{zurich} for more
references).
\begin{cor}
Let $\left(B_{t}\right)$ denote a standard Brownian Motion and let
$$g\equiv sup\left\{t\leq1:\;B_{t}=0\right\}.$$ Then:
$$\mathbb{P}\left[g>t\mid\mathcal{F}_{t}\right]=\sqrt{\dfrac{2}{\pi}}\int_{\frac{\left|B_{t}\right|}{\sqrt{1-t}}}^{\infty}
dy\exp\left(-\frac{y^{2}}{2}\right),$$and
$$A_{t}^{g}=\sqrt{\dfrac{2}{\pi}}\int_{0}^{t\wedge
1}\dfrac{dL_{u}}{\sqrt{1-u}}.$$
\end{cor}
\begin{proof}
It suffices to take $\mu\equiv\dfrac{1}{2}$ in Proposition
\ref{surmartdegamma}.
\end{proof}
It is well known that $g$ is arc sine distributed (see for example
\cite{revuzyor}); with the help of proposition \ref{surmartdegamma},
we can recover this result and extend it to the case of any Bessel
process with dimension $2(1-\mu)$, thus recovering a result also
obtained and proved by Barlow, Pitman and Yor
(\cite{barlowpitmanyor}) using excursion theory (see also Dynkin
\cite{dynkin}):
\begin{cor}
The variable $g_{\mu}$ follows the law:
$$\mathbb{P}\left(g_{\mu}\in
dt\right)=\dfrac{\sin\left(\mu\pi\right)}{\pi}\dfrac{dt}{t^{1-\mu}\left(1-t\right)^{\mu}},\;0<t<1,$$i.e.
the Beta law with parameters $(\mu,1-\mu)$. In particular,
$\mathbb{P}\left(g\in
dt\right)=\dfrac{1}{\pi}\dfrac{dt}{\sqrt{t\left(1-t\right)}},$ i.e.
$g$ is arc sine distributed.
\end{cor}
\begin{proof}
From Proposition \ref{surmartdegamma}, $(2)$, for every Borel
function $f:\;\left[0,1\right]\rightarrow\mathbb{R}_{+}$, we have:
\begin{equation}\label{interne}
\mathbb{E}\left[f\left({g_{\mu}}\right)\right]=\dfrac{1}{2^{\mu}\mu\Gamma\left(\mu\right)}\mathbb{E}\left[\int_{0}^{1}dL_{u}\dfrac{f\left(u\right)}{\left(1-u\right)^{\mu}}\right]=
\dfrac{1}{2^{\mu}\mu\Gamma\left(\mu\right)}
\int_{0}^{1}d_{u}\mathbb{E}\left[L_{u}\right]\dfrac{f\left(u\right)}{\left(1-u\right)^{\mu}}.
\end{equation}
By the scaling property of $\left(L_{t}\right)$,
$$\mathbb{E}\left[L_{u}\right]=u^{\mu}\mathbb{E}\left[L_{1}\right].$$Moreover,
by definition of $\left(L_{t}\right)$,
$$\mathbb{E}\left[L_{1}\right]=\mathbb{E}\left[R_{1}^{2\mu}\right];$$since
$R_{1}^{2}$ is distributed as $2\;gam(1-\mu)$, we have
$$\mathbb{E}\left[R_{1}^{2\mu}\right]=\dfrac{2^{\mu}}{\Gamma\left(1-\mu\right)}.$$Now, plugging this in
(\ref{interne}) yields:
$$\mathbb{E}\left[f\left({g_{\mu}}\right)\right]=\dfrac{1}{\Gamma\left(\mu\right)\Gamma\left(1-\mu\right)}
\int_{0}^{1}du\dfrac{f\left(u\right)}{u^{1-\mu}\left(1-u\right)^{\mu}}.$$To
conclude, it suffices to use the duplication formula for the Gamma
function (\cite{AARoy}):
$$\Gamma\left(\mu\right)\Gamma\left(1-\mu\right)=\dfrac{\pi}{\sin\left(\mu\pi\right)}.$$
\end{proof}
We now state a lemma that we shall often use in the sequel:
\begin{lem}[Az\'{e}ma \cite{azema}]\label{azemgeneral}
Let $L$ be an honest time that avoids $\left(\mathcal{F}_{t}\right)$
stopping times, i.e. for every $\left(\mathcal{F}_{t}\right)$
stopping time $T$, we have $\mathbb{P}\left[L=T\right]=0$, and let
$$Z_{t}^{L}\equiv\mathbb{P}\left[L>t\mid\mathcal{F}_{t}\right].$$ Let $$Z_{t}^{L}=M_{t}^{L}-A_{t},$$ denote its Doob-Meyer decomposition. Then $A_{\infty}$ follows the exponential law with parameter $1$ and
the measure $dA_{t}$ is carried by the set
$\left\{t:\;Z_{t}^{L}=1\right\}$. Moreover, $A$ does not increase
after $L$, i.e. $A_{L}=A_{\infty}$. We also have:
$$L=\sup\left\{t:\;1-Z_{t}^{L}=0\right\}.$$
\end{lem}
\begin{cor}
The variable
$$\frac{1}{2^{\mu}\Gamma\left(1+\mu\right)}\int_{0}^{1}\dfrac{dL_{u}}{\left(1-u\right)^{\mu}}$$is
exponentially distributed with expectation $1$; consequently, its
law is independent of $\mu$.
\end{cor}
\begin{proof}
The random time $g_{\mu}$ is honest by definition (it is the end of
a predictable set). It also avoids stopping times since
$A_{t}^{g_{\mu}}$ is continuous (this can also be seen as a
consequence of the strong Markov property for $R$ and the fact that
$0$ is instantaneously reflecting). Thus the result of the corollary
is a consequence of Proposition \ref{surmartdegamma} and Lemma
\ref{azemgeneral}.
\end{proof}
We can also use Proposition \ref{surmartdegamma} to give an example
of a remarkable random time, called pseudo-stopping time
(\cite{AshkanYor}), which is not a stopping time but which satisfies
the stopping theorem:
\begin{cor}
Define:
$$\rho\equiv\sup\left\{t<g_{\mu}:\;\dfrac{R_{t}}{\sqrt{1-t}}=\sup_{u<g_{\mu}}\dfrac{R_{u}}{\sqrt{1-u}}\right\}.$$
Then, $\rho$ is a pseudo-stopping time, i.e. for every
$\left(\mathcal{F}_{t}\right)$ uniformly integrable martingale
$\left(M_{t}\right)$, we have:
$$\mathbb{E}\left[M_{\rho}\right]=\mathbb{E}\left[M_{\infty}\right].$$
\end{cor}
\begin{proof}
It is an easy consequence of Proposition \ref{surmartdegamma} and
Proposition 5 in \cite{AshkanYor}.
\end{proof}
\section{The conditional law of the last zero before time $1$}
In this section, for simplicity, we consider
$\left(R_{t}\right)_{t\leq1}$; one could easily replace $1$ with any
fixed time $T$. Our aim in this section is twofold:
\begin{itemize}
\item to investigate further for the distribution of $g_{\mu}$ by
computing its conditional distribution given $\mathcal{F}_{t}$;
\item to test the stopping theorem on martingales whose terminal
values are $\sigma\left(g_{\mu}\right)$ measurable, extending thus
the work by Yor for the Brownian Motion in \cite{zurich} (i.e.
$\mu=\dfrac{1}{2}$).
\end{itemize}
We shall need the following lemma which also appears in
\cite{Ashkannonstop} for the resolution of some martingale
equations:
\begin{lem}
Let $L$ be an honest time that avoids stopping times and let
$(K_{t})_{t\geq0}$ be a predictable process such that
$\mathbb{E}\left[\left|K_{L}\right|\right]<\infty$. Then for every
$t\geq0$:
\begin{equation}\label{decoimport}
    \mathbb{E}\left[K_{L}\mid
    \mathcal{F}_{t}\right]=K_{L_{t}}\mathbb{P}\left(L\leq t\mid
    \mathcal{F}_{t}\right)+\mathbb{E}\left[\int_{t}^{\infty}K_{s}dA_{s}\mid
    \mathcal{F}_{t}\right],
\end{equation}where $A_{t}$ is the dual predictable projection of $\mathbf{1}_{\left(L\leq t\right)}$ and $$L_{t}=\sup\left\{s\leq t:\;1-Z_{s}^{L}=0\right\}.$$
Moreover, the latter martingale can also be written as:
$$\mathbb{E}\left[K_{L}\mid
    \mathcal{F}_{t}\right]=-\int_{0}^{t}K_{L_{s}}dM_{s}^{L}+\mathbb{E}\left[\int_{0}^{\infty}K_{s}dA_{s}\mid
    \mathcal{F}_{t}\right],$$where $\left(M_{s}^{L}\right)$ is
    defined in Lemma \ref{azemgeneral}.
\end{lem}
\begin{proof}
\begin{eqnarray*}
  \mathbb{E}\left[K_{L}\mid
    \mathcal{F}_{t}\right] &=& \mathbb{E}\left[K_{L}\mathbf{1}_{L\leq t}\mid
    \mathcal{F}_{t}\right]+\mathbb{E}\left[K_{L}\mathbf{1}_{L> t}\mid
    \mathcal{F}_{t}\right]  \\
   &=& K_{L_{t}}\mathbb{P}\left(L\leq t\mid
    \mathcal{F}_{t}\right)+\mathbb{E}\left[K_{L}\mathbf{1}_{L> t}\mid
    \mathcal{F}_{t}\right].
\end{eqnarray*}Now, let $\Gamma_{t}$ be an
$\left(\mathcal{F}_{t}\right)$ measurable set;
$$\mathbb{E}\left[K_{L}\mathbf{1}_{L>
t}\mathbf{1}_{\Gamma_{t}}\right]=\mathbb{E}\left[\int_{t}^{\infty}K_{s}dA_{s}\mathbf{1}_{\Gamma_{t}}\right];$$hence
$$\mathbb{E}\left[K_{L}\mathbf{1}_{L> t}\mid
    \mathcal{F}_{t}\right]=\mathbb{E}\left[\int_{t}^{\infty}K_{s}dA_{s}\mid
    \mathcal{F}_{t}\right],$$and this completes the proof of the
    first part of the
    lemma. The second part follows from balayage arguments; indeed:
    \begin{eqnarray*}
      K_{L_{t}}\mathbb{P}\left(L\leq t\mid
    \mathcal{F}_{t}\right) &=& K_{L_{t}}\left(1-Z_{t}^{L}\right) \\
       &=& -\int_{0}^{t}K_{L_{s}}dM_{s}^{L}+\int_{0}^{t}K_{s}dA_{s}.
    \end{eqnarray*}Now, since $\mathbb{E}\left[\int_{t}^{\infty}K_{s}dA_{s}\mid
    \mathcal{F}_{t}\right]=\mathbb{E}\left[\int_{0}^{\infty}K_{s}dA_{s}\mid
    \mathcal{F}_{t}\right]-\int_{0}^{t}K_{s}dA_{s}$, we have
    $$K_{L_{t}}\mathbb{P}\left(L\leq t\mid
    \mathcal{F}_{t}\right)+\mathbb{E}\left[\int_{t}^{\infty}K_{s}dA_{s}\mid
    \mathcal{F}_{t}\right]=-\int_{0}^{t}K_{L_{s}}dM_{s}^{L}+\mathbb{E}\left[\int_{0}^{\infty}K_{s}dA_{s}\mid
    \mathcal{F}_{t}\right],$$ and the proof of the lemma is now
    complete.
\end{proof}
Now we shall obtain some closed formulae for martingales whose
terminal values are $g_{\mu}$ measurable and hence obtain the
conditional laws of $g_{\mu}$ given $\mathcal{F}_{t},\;t\geq0$. Let
$h:\;\left[0,1\right]\rightarrow\mathbb{R}_{+}$ be a Borel function
and define:
$$M_{t}^{h}\equiv\mathbb{E}\left[h\left(g_{\mu}\right)\mid\mathcal{F}_{t}\right].$$
The problem of computing the martingales $\left(M_{t}^{h}\right)$
can be dealt with Lemma \ref{decoimport}, which takes here the
following form (recall our processes are stopped at $1$): for any
nonnegative predictable process $(K_{t})_{t\leq 1}$ (such that
$\mathbb{E}\left[K_{g_{\mu}}\right]<\infty$),
\begin{equation}\label{propanulairebis}
    \mathbb{E}\left[K_{g_{\mu}}\mid
    \mathcal{F}_{t}\right]=K_{g_{\mu}\left(t\right)}\mathbb{P}\left(g_{\mu}\leq t\mid
    \mathcal{F}_{t}\right)+\mathbb{E}\left[\int_{t}^{1}K_{s}dA_{s}^{g_{\mu}}\mid
    \mathcal{F}_{t}\right],
\end{equation}
where
$$g_{\mu}\left(t\right)=\sup\left\{s\leq t:\;Z_{s}^{g_{\mu}}=1\right\}=\sup\left\{s\leq t:\;R_{s}=0\right\}.$$
The fact that $\sup\left\{s\leq
t:\;Z_{s}^{g_{\mu}}=1\right\}=\sup\left\{s\leq t:\;R_{s}=0\right\}$
can be seen on the expression of $Z_{t}^{g_{\mu}}$ in Proposition
\ref{surmartdegamma}.

The following lemma will help us to compute explicitly the quantity
$\mathbb{E}\left[\int_{t}^{1}K_{s}dA_{s}^{g_{\mu}}\mid
    \mathcal{F}_{t}\right]$ when $K$ is a deterministic function:
\begin{lem}\label{lemsuivangul}
Let $\left(L_{t}^{a}\right)$ denote the local time at
$a\in\mathbb{R}_{+}$ of the Bessel Process $\left(R_{t}\right)$.
\begin{enumerate}
\item For every Borel function
$h:\mathbb{R}_{+}\rightarrow\mathbb{R}_{+}$,
\begin{equation}\label{katsuc}
    \mathbb{E}\left[\int_{t}^{1}dL_{u}^{a}h\left(u\right)\mid\mathcal{F}_{t}\right]=
    \int_{t}^{1}duh\left(u\right)p^{(\nu)}\left(u-t;R_{t},a\right);
\end{equation}
\item Consequently, for every Borel function
$h:\mathbb{R}_{+}\rightarrow\mathbb{R}_{+}$, we have:
\begin{equation}
\mathbb{E}\left[\int_{t}^{1}dA_{u}^{g_{\mu}}h\left(u\right)\mid\mathcal{F}_{t}\right]=
    \dfrac{\sin\left(\pi\mu\right)}{\pi}\int_{0}^{1}dz\dfrac{h\left(t+z\left(1-t\right)\right)}{\left(1-z\right)^{\mu}z^{1-\mu}}\exp\left(-\dfrac{R_{t}^{2}}{2z\left(1-t\right)}\right).
\end{equation}
\end{enumerate}
\end{lem}
\begin{proof}
$(1)$. First, by the generalized occupation density formula
(\ref{occupationgeneralisee}), for every nonnegative Borel function
$f$, we have:
\begin{equation}\label{rrr}
    \int_{t}^{1}duf\left(R_{u}\right)h\left(u\right)=\int_{0}^{\infty}m\left(da\right)f\left(a\right)\int_{t}^{1}dL_{u}^{a}h\left(u\right).
\end{equation}We also have:
\begin{eqnarray*}
  \mathbb{E}\left[\int_{t}^{1}duf\left(R_{u}\right)h\left(u\right)\mid\mathcal{F}_{t}\right] &=& \int_{t}^{1}duh\left(u\right)\mathbb{E}\left[f\left(R_{u}\right)\mid\mathcal{F}_{t}\right] \\
   &=& \int_{0}^{\infty}m\left(da\right)f\left(a\right)\int_{t}^{1}duh\left(u\right)p^{(\nu)}\left(u-t;R_{t},a\right)
\end{eqnarray*}
and from (\ref{rrr}) we obtain:
$$\mathbb{E}\left[\int_{t}^{1}dL_{u}^{a}h\left(u\right)\mid\mathcal{F}_{t}\right]=
    \int_{t}^{1}duh\left(u\right)p^{(\nu)}\left(u-t;R_{t},a\right).$$
$(2)$. We know from Proposition (\ref{surmartdegamma}) that:
$A_{t}^{g_{\mu}}=\frac{1}{2^{\mu}\mu\Gamma\left(\mu\right)}\int_{0}^{t\wedge
1}\dfrac{dL_{u}}{\left(1-u\right)^{\mu}}$. Plugging this into
(\ref{katsuc}) yields ($a=0$):
$$\mathbb{E}\left[\int_{t}^{1}dA_{u}^{g_{\mu}}h\left(u\right)\mid\mathcal{F}_{t}\right]=\int_{t}^{1}duh\left(u\right)
\dfrac{\exp\left(-\frac{R_{t}^{2}}{2\left(u-t\right)}\right)}{\Gamma\left(\mu\right)
\Gamma\left(1-\mu\right)\left(1-u\right)^{\mu}\left(t-u\right)^{1-\mu}}.$$Now,
using the duplication formula: $$\Gamma\left(\mu\right)
\Gamma\left(1-\mu\right)=\dfrac{\pi}{\sin\left(\mu\pi\right)},$$and
making the change of variable $u=t+z\left(1-t\right)$, we obtain:
$$\mathbb{E}\left[\int_{t}^{1}dA_{u}^{g_{\mu}}h\left(u\right)\mid\mathcal{F}_{t}\right]=
    \dfrac{\sin\left(\pi\mu\right)}{\pi}\int_{0}^{1}dz\dfrac{h\left(t+z\left(1-t\right)\right)}{\left(1-z\right)^{\mu}z^{1-\mu}}\exp\left(-\dfrac{R_{t}^{2}}{2z\left(1-t\right)}\right),$$
    which completes the proof of the lemma.
\end{proof}
Now, we shall give three nice corollaries of Lemma
\ref{lemsuivangul}. The last two corollaries  extend naturally some
results of Yor (\cite{columbia}) for the Brownian Motion to any
Bessel process of dimension $\delta \in (0,2)$.
\begin{cor}\label{decodemartarrete}
Let $h:[0,1]\rightarrow\mathbb{R}_{+}$, be a Borel function, then:
\begin{equation*}
    \mathbb{E}\left[h\left(g_{\mu}\right)\mid
    \mathcal{F}_{t}\right]=h\left(g_{\mu}\left(t\right)\right)\left(1-Z_{t}^{g_{\mu}}\right)+\mathbb{E}\left[h\left(g_{\mu}\right)\mathbf{1}_{\left(g_{\mu}>t\right)}\mid \mathcal{F}_{t}\right];
\end{equation*}with
\begin{equation}\label{equatioonpont}
    \mathbb{E}\left[h\left(g_{\mu}\right)\mathbf{1}_{\left(g_{\mu}>t\right)}\mid \mathcal{F}_{t}\right]=
   \dfrac{\sin\left(\pi\mu\right)}{\pi}\int_{0}^{1}dz\dfrac{h\left(t+z\left(1-t\right)\right)}{\left(1-z\right)^{\mu}z^{1-\mu}}\exp\left(-\dfrac{R_{t}^{2}}{2z\left(1-t\right)}\right).
\end{equation}Consequently, the law of $g_{\mu}$ given
$\mathcal{F}_{t}$, which we note $\lambda_{t}\left(dz\right)$, is
given by:
$$\lambda_{t}\left(dz\right)=\left(1-Z_{t}^{g_{\mu}}\right)\varepsilon_{g_{\mu}\left(t\right)}\left(dz\right)+\mathbf{1}_{\left(t,1\right)}\left(z\right)\dfrac{\sin\left(\pi\mu\right)}{\pi}\dfrac{\exp\left(-\frac{R_{t}^{2}}{2\left(z-t\right)}\right)}{\left(1-z\right)^{\mu}\left(z-t\right)^{1-\mu}}dz,$$
and taking $t=0$, we recover the generalized arc sine law.
\end{cor}
As a consequence of this corollary, we can also see how the stopping
theorem fails to hold for the family of martingales $M_{t}^{h}\equiv
\mathbb{E}\left[h\left(g_{\mu}\right)\mid \mathcal{F}_{t}\right]$,
thus completing the examples in \cite{zurich} and
\cite{Ashkannonstop}:
\begin{cor}\label{corporgam}
Let $h:[0,1]\rightarrow\mathbb{R}_{+}$, be a Borel function, and
define $M_{t}^{h}=\mathbb{E}\left[h\left(g_{\mu}\right)\mid
\mathcal{F}_{t}\right]$; then
\begin{equation}
    \mathbb{E}\left[M_{\infty}^{h}\mid
\mathcal{F}_{g_{\mu}}\right]=h\left(g_{\mu}\right),
\end{equation}whilst
\begin{equation}
    M_{g_{\mu}}^{h}=\dfrac{\sin\left(\pi\mu\right)}{\pi}\int_{0}^{1}dz\dfrac{1}{\left(1-z\right)^{\mu}z^{1-\mu}}h\left(g_{\mu}+z\left(1-g_{\mu}\right)\right).
\end{equation}
\end{cor}
We can also compute explicitly the martingale
$\mathbb{E}\left[L_{1}\mid\mathcal{F}_{t}\right]$:
\begin{cor}
Let $$X_{t}\equiv\mathbb{E}\left[L_{1}\mid\mathcal{F}_{t}\right].$$
Then,
$$X_{1}\equiv X_{\infty}=L_{1},$$whilst$$X_{g_{\mu}}=L_{1}+\dfrac{2^{\mu}}{\Gamma\left(1-\mu\right)}\left(1-g_{\mu}\right)^{\mu}.$$
\end{cor}
\begin{proof}
It suffices to take $h\equiv 1$ in (\ref{katsuc}).
\end{proof}
Taking $\mu=\dfrac{1}{2}$, we recover the following results obtained
for  Brownian Motion in \cite{columbia}:
\begin{cor}
Let $\left(B_{t}\right)$ be the standard Brownian Motion, and denote
$\left(\ell_{t}\right)$ its local time at zero. Let $g\equiv
g_{1/2}$. Then, for any Borel function
$h:[0,1]\rightarrow\mathbb{R}_{+}$, the following identities hold:
\begin{enumerate}
\item
$$\left.\mathbb{E}\left[h\left(g\right)\mid\mathcal{F}_{t}\right]\right|_{t=g}=\dfrac{1}{\pi}\int_{0}^{1}\dfrac{dz}{\sqrt{z\left(1-z\right)}}h\left(g+z\left(1-g\right)\right).$$
\item
$$\left.\mathbb{E}\left[\ell_{1}\mid\mathcal{F}_{t}\right]\right|_{t=g}=\ell_{1}+\sqrt{\dfrac{2}{\pi}}\sqrt{1-g}.$$
\end{enumerate}
\end{cor}
\section{Some martingales with the same set of zeros as $R_{t}$ and which satisfy the stopping theorem with respect to $g_{\mu}$}
In this section, we shall illustrate with some examples related to
Bessel processes the seminal work of Az\'{e}ma and Yor
(\cite{azemayorzero}) on zeros of continuous martingales. Yor has
specialized further this work to the important case of the standard
Brownian Motion, giving explicit examples of martingales which have
the same zeros as $\left(B_{t}\right)_{t\leq1}$ (\cite{zurich},
chapter 14). Quite unexpectedly, the study of martingales which have
 same zeros leads to some discussion on the stopping theorem when
stopping times are replaced with honest times
(\cite{azemajeulinknightyor}, \cite{Ashkannonstop}). Our aim here is
to provide a one parameter ($\mu$) extension of some results in
\cite{zurich}, chapter 14, and to give more examples of martingales
associated with honest times satisfying the stopping theorem (which
is, regarding the results of the previous section, quite
exceptional). Even though $\left(R_{t}\right)$ is not a martingale
(and not even  a semimartingale when $\mu>\dfrac{1}{2}$), the
methods of Az\'{e}ma and Yor apply remarkably here. More precisely:
\begin{prop}\label{zerosdebeesel}
Define:
$$\mathcal{Z}_{1}=\left\{t\in[0,1]:\;R_{t}=0\right\},$$
and let $\left(M_{t}\right)$ be an $\left(\mathcal{F}_{t}\right)$
martingale. Then the following are equivalent:
\begin{enumerate}
\item for all $t\in \mathcal{Z}_{1}$,\;$M_{t}=0$;
\item $M_{g_{\mu}}=0$.

\noindent If furthermore $\left(M_{t}\right)$ is uniformly
integrable, then the previous assertions are equivalent to:
\item
$\mathbb{E}\left[M_{\infty}\mid\mathcal{F}_{g_{\mu}}\right]=0$;
consequently, in this case, we have:
\begin{equation}\label{knightmaison}
\mathbb{E}\left[M_{\infty}\mid\mathcal{F}_{g_{\mu}}\right]=M_{g_{\mu}}.
\end{equation}
\end{enumerate}
\end{prop}
\begin{proof}
We only prove $(1)\Leftrightarrow (2)$, as the proof of the
equivalence with $(3)$ can be found in
\cite{azemayorzero,azemajeulinknightyor,Ashkannonstop}.

Only the implication $(2)\Rightarrow (1)$ is not obvious. Assume
that $M_{g_{\mu}}=0$; then, we also have
$\left|M_{g_{\mu}}\right|=0$. But from Proposition
\ref{surmartdegamma},
$$0=\mathbb{E}\left[\left|M_{g_{\mu}}\right|\right]=\mathbb{E}\left[\int_{0}^{1}\left|M_{u}\right|dA_{u}^{g_{\mu}}\right],$$
and consequently
$$\left|M_{u}\right|=0,\;dA_{u}^{g_{\mu}}d\mathbb{P}\;a.s.$$ and from
Lemma \ref{azemgeneral} and Proposition \ref{surmartdegamma}, this
means that $M_{t}=0$ if and only if $t\in \mathcal{Z}_{1}$.
\end{proof}
Now, as a consequence of Proposition \ref{zerosdebeesel}, we can
associate canonically with a uniformly integrable martingale another
uniformly integrable martingale whose zeros are in $\mathcal{Z}_{1}$
and which satisfy (\ref{knightmaison}). More precisely, if
$\left(M_{t}\right)$ is a uniformly integrable martingale, then the
uniformly integrable martingale $\left(\widehat{M}_{t}\right)$
defined by:
$$\widehat{M}_{t}\equiv\mathbb{E}\left[\left(M_{\infty}-\mathbb{E}\left[M_{\infty}\mid\mathcal{F}_{g_{\mu}}\right]\right)\mid\mathcal{F}_{t}\right],$$
satisfies the equivalent assertions of Proposition
\ref{zerosdebeesel} ($(3)$). We shall now illustrate this fact on an
example reminiscent of Yor's study for the standard Brownian Motion
(\cite{zurich,columbia}). More precisely, with
$f:\mathbb{R}_{+}\rightarrow\mathbb{R}_{+}$, a Borel function, we
associate the martingale:
$$M_{t}^{f}\equiv\mathbb{E}\left[f\left(R_{1}\right)\mid\mathcal{F}_{t}\right].$$
To give an explicit expression for the associated martingale
$\widehat{M}_{t}^{f}$, we shall need some results of Yor
 about generalized meanders as exposed in \cite{zurichI}.
 For a fixed time $T>0$, we call:
$$m_{\mu}^{T}\left(u\right)=\dfrac{1}{\sqrt{T-g_{\mu}\left(T\right)}}R_{g_{\mu}\left(T\right)+u\left(T-g_{\mu}\left(T\right)\right)},\;u\leq1$$the
Bessel meander associated to the Bessel process $R$ of dimension
$2(1-\mu);\;\mu\in(0,1)$.
$\left(m_{\mu}^{T}\left(u\right)\right)_{u\leq 1}$ is independent of
$\mathcal{F}_{g_{\mu}\left(T\right)}$, and the law of
$m_{\mu}^{T}\equiv m_{\mu}^{T}\left(1\right)$ does not depend on
$\mu$ and $T$, and is distributed as the two dimensional Bessel
process at time $1$, i.e.:
\begin{equation}\label{loimeandre}
    \mathbb{P}\left(m_{\mu}^{T}\in
    dx\right)=x\exp\left(-\dfrac{x^{2}}{2}\right)dx.
\end{equation}In the sequel, we shall simply note $m_{\mu}$ for
$m_{\mu}^{1}\left(1\right)$.

Now, let us come back to the case of the martingales
$M_{t}^{f}=\mathbb{E}\left[f\left(R_{1}\right)\mid\mathcal{F}_{t}\right].$
From the Markov property, we have:
$$M_{t}^{f}=\int_{0}^{\infty}m\left(dz\right)f\left(z\right)p^{(\nu)}\left(1-t;R_{t},z\right),$$and
from the properties recalled above about Bessel meanders, we have:
\begin{eqnarray*}
  \mathbb{E}\left[f\left(R_{1}\right)\mid\mathcal{F}_{g_{\mu}}\right] &=& \mathbb{E}\left[f\left(m_{\mu}\sqrt{1-g_{\mu}}\right)\mid\mathcal{F}_{g_{\mu}}\right] \\
   &=&
   \int_{0}^{\infty}dzf\left(z\sqrt{1-g_{\mu}}\right)z\exp\left(-\dfrac{z^{2}}{2}\right).
\end{eqnarray*}Now, with the help of Corollary \ref{decodemartarrete}, we
are able to compute
$\mathbb{E}\left[\mathbb{E}\left[f\left(R_{1}\right)\mid\mathcal{F}_{g_{\mu}}\right]\mid\mathcal{F}_{t}\right]$:
\begin{prop}\label{longue}
Define as above:
$$\widehat{M}_{t}^{f}\equiv
\mathbb{E}\left[\left(f\left(R_{1}\right)-\mathbb{E}\left[f\left(R_{1}\right)\mid\mathcal{F}_{g_{\mu}}\right]\right)\mid\mathcal{F}_{t}\right].$$
Define:
\begin{eqnarray*}
  g_{\mu}\left(t\right) &\equiv& \sup\left\{s\leq t:\;R_{s}=0\right\}; \\
  \theta_{\mu}\left(x\right) &=&
  \dfrac{1}{2^{\mu-1}\Gamma\left(\mu\right)}\int_{x}^{\infty}dzz^{2\mu-1}\exp\left(-\dfrac{z^{2}}{2}\right).
\end{eqnarray*}
Then $\left(\widehat{M}_{t}^{f}\right)$ satisfies the equivalent
assertions of Proposition \ref{zerosdebeesel} and can be expressed
as:
$$\widehat{M}_{t}^{f}=\widehat{M}_{t}^{(1)}-\widehat{M}_{t}^{(2)}-\widehat{M}_{t}^{(3)},$$where:
\begin{eqnarray*}
  \widehat{M}_{t}^{(1)} &=& \int_{0}^{\infty}m\left(dz\right)f\left(z\right)p^{(\nu)}\left(1-t;R_{t},z\right), \\
  \widehat{M}_{t}^{(2)} &=& \theta_{\mu}\left(\dfrac{R_{t}}{\sqrt{1-t}}\right)\int_{0}^{\infty}dzzf\left(z\sqrt{1-g_{\mu}\left(t\right)}\right)\exp\left(-\dfrac{z^{2}}{2}\right), \\
  \widehat{M}_{t}^{(3)} &=& \dfrac{\sin\left(\pi\mu\right)}{\pi}\int_{0}^{\infty}dzz\exp\left(-\dfrac{z^{2}}{2}\right)
  \int_{0}^{1}dw\dfrac{f\left(z\sqrt{1-t}\sqrt{1-w}\right)}{\left(1-w\right)^{\mu}w^{\mu}}
  \exp\left(-\dfrac{R_{t}^{2}}{2w\left(1-t\right)}\right).
\end{eqnarray*}
\end{prop}
\begin{proof}
The proof is obtained with the help of Corollary
\ref{decodemartarrete}, Proposition \ref{surmartdegamma} and a few
elementary computations.
\end{proof}
\begin{rem}
The above proposition tells us that though Proposition
\ref{zerosdebeesel} is simple, obtaining the explicit expression of
the projections on $\mathcal{F}_{g_{\mu}}$ and then on
$\mathcal{F}_{t}$ can be difficult in practice.
\end{rem}\bigskip
Now, we shall use the remarkable fact that the density of the
generalized meander at time $1$ does not depend on $\mu$ and
satisfies a functional equation to build a another family of
martingales with the same properties as the martingales
$\left(\widehat{M}_{t}^{f}\right)$.
\begin{prop}
Let $f:\mathbb{R}_{+}\rightarrow\mathbb{R}$, be a function of class
$\mathcal{C}^{2}$, with compact support. Define
$$X^{f}\equiv
f\left(R_{1}\right)-f\left(0\right)-R_{1}f'\left(R_{1}\right)+\left(1-g_{\mu}\right)f''\left(R_{1}\right),$$
and
$$X_{t}^{f}\equiv \mathbb{E}\left[X^{f}\mid
\mathcal{F}_{t}\right].$$Then, $\left(X_{t}^{f}\right)$ satisfies
(\ref{knightmaison}), i.e. $$\mathbb{E}\left[X_{\infty}^{f}\mid
\mathcal{F}_{g_{\mu}}\right]=X_{g_{\mu}}^{f},$$and
$$X_{t}^{f}=0\Leftrightarrow t\in\mathcal{Z}_{1}.$$
\end{prop}
\begin{proof}
As recalled at the beginning of the previous subsection,
$\mathbb{P}\left(m_{\mu}\in
d\rho\right)=\rho\exp\left(-\dfrac{\rho^{2}}{2}\right)d\rho$. It is
not difficult to show that for
$f:\mathbb{R}_{+}\rightarrow\mathbb{R}$, a function of class
$\mathcal{C}^{2}$, which is compactly supported, we have:
$$\mathbb{E}\left[f\left(m_{\mu}\right)\right]=\mathbb{E}\left[f\left(0\right)+m_{\mu}f'\left(m_{\mu}\right)-f''\left(m_{\mu}\right)\right].$$Now,
replacing $f$ with $f\left(\kappa\bullet\right)$ in the above, with
$\kappa\in\mathbb{R}$, we obtain:
$$\mathbb{E}\left[f\left(\kappa m\right)\right]=\mathbb{E}\left[f\left(0\right)+\kappa mf'\left(\kappa m\right)-\kappa^{2}f''\left(\kappa m\right)\right].$$
Next, as $R_{1}=\sqrt{1-g_{\mu}}m_{\mu}$, with $m_{\mu}$ independent
from $\mathcal{F}_{g_{\mu}}$, we have:
$$\mathbb{E}\left[X\mid\mathcal{F}_{g_{\mu}}\right]=0,$$and the result of the Proposition follows from Proposition
\ref{zerosdebeesel}.
\end{proof}
\begin{rem}
In fact, the result of the Proposition is still true if we only
assume that:
$\mathbb{E}\left[\left|f\left(m_{\mu}\right)\right|\right]<\infty,\;\mathbb{E}\left[m_{\mu}\left|f'\left(m_{\mu}\right)\right|\right]<\infty,\;\mathbb{E}\left[\left|f''\left(m_{\mu}\right)\right|\right]<\infty$.
\end{rem}
\begin{rem}
More generally, if $\varphi:\mathbb{R}_{+}\times
[0,1]\rightarrow\mathbb{R}_{+}$, then, using the properties of the
meander, we obtain that:
$$\mathbb{E}\left[\varphi\left(R_{1},g_{\mu}\right)\mid
\mathcal{F}_{g_{\mu}}\right]=0,$$if and only if:
$$\int_{0}^{\infty}dz\;z\exp\left(-\dfrac{z^{2}}{2\left(1-g_{\mu}\right)}\right)\varphi\left(z,g_{\mu}\right)=0.$$
\end{rem}
\section{A one parameter extension of Az\'{e}ma's second martingale}
In this section, we shall associate with the pair
$\left(R_{t},L_{t}\right)$ a family of local martingales reminiscent
of the Az\'{e}ma-Yor martingales and use them to compute the
distribution of the local time at some stopping times. Furthermore,
we shall project these martingales on the filtration of the zeros of
$\left(R_{t}\right)$ to obtain some remarkable martingales; in
particular, we prove a one parameter extension of Az\'{e}ma's second
martingale, i.e. we find a submartingale in the filtration of the
zeros, which has the same local time $\left(L_{t}\right)$ as
$\left(R_{t}\right)$.

Let $f:\mathbb{R}_{+}\rightarrow\mathbb{R}$ be a locally bounded
Borel function and let
$F\left(x\right)=\int_{0}^{x}dzf\left(z\right)$.
\begin{prop}\label{martreflechies}
Let $f$ and $F$ be defined as above. Then the process
$$F\left(L_{t}\right)-f\left(L_{t}\right)R_{t}^{2\mu},$$is a local
martingale.
\end{prop}
\begin{proof}
We have:
$$R_{t}^{2\mu}=N_{t}+L_{t},$$and $dL_{t}$ is carried by
$\left\{t:\;R_{t}=0\right\}$.
Hence as an application of Skorokhod's reflection lemma (see \cite{revuzyor}, chapter VI) yields:%
\begin{equation*}
L_{t}=\sup_{s\leq t}\left( -N_{s}\right)
\end{equation*}
Then applying the balayage formula (see \cite{revuzyor}, chapter VI,
p.262) we obtain:
\begin{equation}\label{mlocbalay}
F\left( L_{t}\right) -f\left( L_{t}\right)R_{t}^{2\mu}
\end{equation}%
is a local martingale.
\end{proof}
Now, for $a>0$, let:
$$T_{a}\equiv\inf\left\{t:\;R_{t}=a\right\}.$$
\begin{cor}
The random variable $L_{T_{a}}$ is distributed as an exponential
variable of parameter $a^{2\mu}$.
\end{cor}
\begin{proof}
Taking $F\left(x\right)\equiv\exp\left(-\theta x\right),\;\theta>0$
in Proposition \ref{martreflechies}, we obtain that:
$$M_{t}\equiv\exp\left(-\theta L_{t\wedge T_{a}}\right)\left(1+\theta R_{t\wedge T_{a}}^{2\mu}\right),$$
is a local martingale. Since it is bounded, it is a uniformly
integrable martingale and the optional stopping theorem yields:
$$\mathbb{E}\left[\exp\left(-\theta L_{T_{a}}\right)\right]=\dfrac{1}{1+\theta
a^{2\mu}},$$and thus $L_{T_{a}}$ is distributed as an exponential
variable of parameter $a^{2\mu}$.
\end{proof}
\begin{rem}
When $\mu=\dfrac{1}{2}$, we recover a well known result for the
standard Brownian Motion: in this special case,
$T_{a}\equiv\inf\left\{t:\;\left|B_{t}\right|=a\right\}.$
\end{rem}\bigskip
Now, we shall project the martingales of Proposition
\ref{martreflechies} on the smaller filtration
$\left(\mathcal{F}_{g_{\mu}\left(t\right)}\right)$.
\begin{lem}
For $u\leq v$, we have:
$$\mathcal{F}_{g_{\mu}\left(u\right)}\subset\mathcal{F}_{g_{\mu}\left(v\right)},$$
and for every $t\geq0$, we have:
$$\mathcal{F}_{g_{\mu}\left(t\right)}\subset\mathcal{F}_{t}.$$
\end{lem}
\begin{proof}
This is a consequence of Theorem 27, p. 142 in \cite{delmaismey}
about the sigma algebras associated with honest times.
\end{proof}
Our aim now is to find martingales for the filtration
$\left(\mathcal{F}_{g_{\mu}\left(t\right)}\right)$, the filtration
of the zeros of $R$. There exist some general results, essentially
due to Az\'{e}ma, on excursion theory for closed optional random
sets (\cite{azemaferme}, \cite{delmaismey}), and our framework would
fit in that general setting. However, the computations are not
always easy to carry in the latter setting and we shall use more
elementary arguments to deal with this problem. Az\'{e}ma did some
exact computations in the Brownian setting and obtained the
celebrated Az\'{e}ma's martingales. We shall obtain here a one
parameter generalization of Az\'{e}ma's second martingale.

We are interested in local martingales of Proposition
\ref{martreflechies} which are true martingales. This happens for
example if $F$ and $F'=f$ are with compact support, or if $f$ is a
probability density on $\mathbb{R}_{+}$.
\begin{prop}\label{generazemainfini}
Let $f$ be chosen such that the local martingale
$M_{t}=f\left(L_{t}\right)R_{t}^{2\mu}-F\left(L_{t}\right)$ is a
true martingale. Define:
$$\Lambda_{t}\equiv\mathbb{E}\left[M_{t}\mid\mathcal{F}_{g_{\mu}\left(t\right)}\right].$$
Then $\left(\Lambda_{t}\right)$ is a martingale in the filtration
$\left(\mathcal{F}_{g_{\mu}\left(t\right)}\right)$ and we have:
$$\Lambda_{t}=2^{\mu}\Gamma\left(1+\mu\right)f\left(L_{t}\right)\left(t-g_{\mu}\left(t\right)\right)^{\mu}-F\left(L_{t}\right).$$
\end{prop}
\begin{proof}
First, we note that $L_{g_{\mu}\left(t\right)}=L_{t}$ ($L$ increases
on the set of zeros of $R$). Consequently, we have:
$$\mathbb{E}\left[M_{t}\mid\mathcal{F}_{g_{\mu}\left(t\right)}\right]=f\left(L_{t}\right)\mathbb{E}\left[R_{t}^{2\mu}\mid\mathcal{F}_{g_{\mu}\left(t\right)}\right]-F\left(L_{t}\right).$$
Now, from the properties of the generalized meander, we have:
$$\mathbb{E}\left[R_{t}^{2\mu}\mid\mathcal{F}_{g_{\mu}\left(t\right)}\right]=\left(t-g_{\mu}\left(t\right)\right)^{\mu}\mathbb{E}\left[\left(m_{\mu}^{(t)}\right)^{2\mu}\right].$$
To conclude, it suffice to note that:
$$\mathbb{E}\left[\left(m_{\mu}^{(t)}\right)^{2\mu}\right]=\int_{0}^{\infty}dz\;z^{2\mu+1}\exp\left(-\frac{z^{2}}{2}\right)=2^{\mu}\Gamma\left(1+\mu\right).$$
\end{proof}
A nice consequence of Proposition \ref{generazemainfini} is the
existence of a remarkable submartingale in the filtration
$\left(\mathcal{F}_{g_{\mu}\left(t\right)}\right)$, which vanishes
on $\mathcal{Z}_{1}$ and which has the same local time as
$\left(R_{t}\right)$. The existence of such a submartingale leads to
the generalization of the so called Az\'{e}ma's second martingale:
by analogy with the Brownian case (see \cite{azemaferme}, p.462), we
shall project the martingale $\left(R_{t}^{2\mu}-L_{t}\right)$ on
the filtration $\left(\mathcal{F}_{g_{\mu}\left(t\right)}\right)$.
\begin{cor}\label{azemunidim}
The stochastic process $$\nu_{t}\equiv
\left(t-g_{\mu}\left(t\right)\right)^{\mu}-c_{\mu}L_{t},$$ where
$$c_{\mu}=\dfrac{1}{2^{\mu}\Gamma\left(1+\mu\right)},$$ is an $\left(\mathcal{F}_{g_{\mu}\left(t\right)}\right)$
martingale. Consequently $Y_{t}\equiv
2^{\mu}\Gamma\left(1+\mu\right)\left(t-g_{\mu}\left(t\right)\right)^{\mu}$
is an $\left(\mathcal{F}_{g_{\mu}\left(t\right)}\right)$
submartingale which vanishes on $\mathcal{Z}_{1}$ and whose local
time is $\left(L_{t}\right)$. In Az\'{e}ma's terminology,
$\left(Y_{t}\right)$ is called the equilibrium submartingale of
$\mathcal{Z}_{1}$, and taking $\mu=\dfrac{1}{2}$, we recover the
calculations of Az\'{e}ma in the Brownian setting (we drop the index
$\mu$ for notational convenience):
\begin{eqnarray*}
  Y_{t} &=& \sqrt{\dfrac{\pi}{2}}\sqrt{t-g_{t}} \\
  \nu_{t} &=& \sqrt{t-g_{t}}-\sqrt{\dfrac{2}{\pi}}\ell_{t},
\end{eqnarray*}the latter martingale being the celebrated
Az\'{e}ma's second martingale (see \cite{zurich}).
\end{cor}
\begin{proof}
It suffices to take $f\equiv1$ in Proposition
\ref{generazemainfini}.
\end{proof}
\begin{rem}
C. Rainer (\cite{rainer}) has some projection formulae for real
diffusions on the slow filtration; our result is not contained in
her work but we could use her results, based on excursion theory, to
extend some of our results to some diffusions on natural scale.
\end{rem}
\begin{cor}
Let $n_{\mu}$ be the L\'{e}vy measure of the inverse local time of
$R$ or the distribution of the lifetime of excursions under the
It\^{o} measure. Then, we have:
$$n_{\mu}\left(dx\right)=\dfrac{1}{2^{\mu}\Gamma\left(\mu\right)}\dfrac{dx}{x^{1+\mu}}.$$
\end{cor}
\begin{proof}
This is a consequence of Corollary \ref{azemunidim} and the
following result of Az\'{e}ma (\cite{azemaferme}, p.452) which takes
in our setting the following form: the stochastic process
$$\overline{M}_{t}\equiv\dfrac{1}{n_{\mu}\left(\left[t-g_{\mu}\left(t\right)\right[\right)}-L_{t}\equiv2^{\mu}\Gamma\left(\mu+1\right)\left(t-g_{\mu}\left(t\right)\right)^{\mu}-L_{t},$$
is a local martingale for the filtration
$\left(\mathcal{F}_{g_{\mu}\left(t\right)}\right)$ (in fact we have
proved it is a true martingale).
\end{proof}\bigskip
Now, to conclude, we shall give a local time estimate for $R$. What
follows is reminiscent of some studies by Knight
(\cite{Knight73,knight78}), Shi (\cite{shi}) and Khoshnevisan
(\cite{Davar}), and is inspired by some recent lectures given by
Marc Yor at Columbia University (\cite{columbia}). It should be
mentioned that our results, which generalize some similar results in
the Brownian setting, can be in fact generalized to a much wider
class of stochastic processes studied in a forthcoming paper
\cite{Ashkanrem}, without assuming any Markov nor scaling property.
This is in fact possible thanks to the martingale techniques we
shall use. More precisely, we shall need the following result,
called Doob's maximal identity. We only mention it without proving
it; the reader can refer to \cite{ashyordoob} for a proof (which is
essentially an application of Doob's optional stopping theorem) and
for some nice applications to enlargements of filtrations and path
decompositions for some large classes of diffusion processes.
\begin{lem}[Doob's maximal
identity] \label{maxeq} Let $\left(N_{t}\right)$ be a continuous and
positive local martingale which satisfies:
$$N_{0}=x,\;x>0;\;\lim_{t\rightarrow\infty}N_{t}=0.$$
If we note $$S_{t}\equiv\sup_{u\leq t}N_{u},$$then, for any $a>0$, we have:%
\begin{enumerate}
\item
\begin{equation} \mathbf{P}\left( S_{\infty }>a\right) =\left(
\frac{x}{a}\right) \wedge 1. \label{loimax}
\end{equation}Hence, $\dfrac{1}{S_{\infty }}$ is a uniform random variable on $%
\left(0,1/a\right)$.
\item For any stopping time $T$:%
\begin{equation}
\mathbf{P}\left( S^{T}>a\mid \mathcal{F}_{T}\right) =\left( \frac{N_{T}}{a}%
\right) \wedge 1 ,  \label{loimaxcond}
\end{equation}%
where
\begin{equation*}
S^{T}=\sup_{u\geq T}N_{u}.
\end{equation*}%
Hence $\dfrac{N_{T}}{S^{T}}$ is also a uniform random variable on
$\left(0,1\right)$, independent of $\mathcal{F}_{T}$.
\end{enumerate}
\end{lem}
Now, we can state and prove our result about local time estimates:
\begin{prop}\label{estimloc}
Let $R$ be a Bessel process of dimension $2(1-\mu)$, with
$\mu\in(0,1)$, $L$ its local time at $0$ (as defined in Section 2).
Define $\tau$ the right continuous inverse of $L$:
$$\tau_{u}=\inf\left\{t\geq0;\;L_{t}> u\right\}.$$
Then, for any $u>0$, and any positive Borel function $\varphi$, we
have:
$$\mathbb{P}\left(\exists t\leq
\tau_{u},\;R_{t}>\varphi\left(L_{t}\right)\right)=1-\exp\left(-\int_{0}^{u}\frac{dx}{\varphi^{2\mu}\left(x\right)}\right),$$and
consequently:
$$\mathbb{P}\left(\exists t\geq
0,\;R_{t}>\varphi\left(L_{t}\right)\right)=1-\exp\left(-\int_{0}^{\infty}\frac{dx}{\varphi^{2\mu}\left(x\right)}\right).$$
\end{prop}
\begin{proof}
For $u>0$, define:
$$\varphi_{u}\left(x\right)=\varphi\left(x\right),\mathrm{if}\;x<u;\;\mathrm{and}\;\varphi_{u}\left(x\right)=\infty\;\mathrm{otherwise}.$$
Now, it is clear that: $$\mathbb{P}\left(\exists t\leq
\tau_{u},\;R_{t}>\varphi\left(L_{t}\right)\right)=\mathbb{P}\left(\exists
t\geq 0,\;R_{t}>\varphi_{u}\left(L_{t}\right)\right).$$ Consider now
the local martingale:
$$M_{t}\equiv F\left(L_{t}\right)-f\left(L_{t}\right)R_{t}^{2\mu},$$where we
choose (with the convention $\frac{1}{\infty}=0$):
$$F\left(x\right)\equiv1-\exp\left(-\int_{x}^{\infty}\frac{dz}{\varphi_{u}^{2\mu}\left(z\right)}\right),$$and
$f=F'$, the Lebesgue derivative of $F$. Now, it is easily checked
that $M$ is a positive local martingale. Moreover,
$\lim_{t\rightarrow\infty}M_{t}=0$. Indeed, since $M$ is a positive
local martingale, it converges almost surely to a limit
$M_{\infty}$. To see that in fact $M_{\infty}=0$, we look at
$\lim_{v\rightarrow\infty}M_{\tau_{v}}$. Since $L_{\tau_{v}}=v$ and
$R_{\tau_{v}}=0$, we easily find that:
$\lim_{v\rightarrow\infty}M_{\tau_{v}}=0$, and consequently,
$M_{\infty}=0$.

Now let us note that if for a given $t_{0}<\infty$, we have
$R_{t_{0}}>\varphi_{u}\left(L_{t_{0}}\right)$, then we must have:
$$M_{t_{0}}>F\left(
L_{t_{0}}\right) -f\left( L_{t_{0}}\right)
\varphi_{u}^{2\mu}\left(L_{t_{0}}\right)=1,$$and hence we easily
deduce from this that:
\begin{eqnarray*}
  \mathbb{P}\left(\exists t\geq
0,\;R_{t}>\varphi_{u}\left(L_{t}\right)\right) &=&
\mathbb{P}\left(\exists t\geq
0,\;R_{t}^{2\mu}>\varphi_{u}^{2\mu}\left(L_{t}\right)\right) \\
                                               &=& \mathbb{P}\left(\sup_{t\geq0}M_{t}>1\right) \\
   &=&  \mathbb{P}\left(\sup_{t\geq0}\dfrac{M_{t}}{M_{0}}>\dfrac{1}{M_{0}}\right)\\
   &=& M_{0},
\end{eqnarray*}where the last equality is obtained by an application
of Doob's maximal identity (Lemma \ref{maxeq}). To conclude, it
suffices to note that
$$M_{0}=1-\exp\left(-\int_{0}^{\infty}\frac{dx}{\varphi_{u}^{2\mu}\left(x\right)}\right)=1-\exp\left(-\int_{0}^{u}\frac{dx}{\varphi_{u}^{2\mu}\left(x\right)}\right).$$
\end{proof}
\begin{cor}
With the hypotheses of the above proposition, the following integral
criterion holds: if
$\int_{0}^{\infty}\frac{dx}{\varphi^{2\mu}\left(x\right)}=\infty$,
then:
\begin{equation*}
\mathbb{P}\left(\forall A>0,\; \exists t\geq
A,\;R_{t}>\varphi\left(L_{t}\right)\right)=1;
\end{equation*}if
$\int_{0}^{\infty}\frac{dx}{\varphi^{2\mu}\left(x\right)}<\infty$,
then: $$\mathbb{P}\left(\forall A>0,\;\exists t\geq
A,\;R_{t}>\varphi\left(L_{t}\right)\right)=0.$$
\end{cor}
\begin{proof}
The event $\left\{\forall A>0,\; \exists t\geq
A,\;R_{t}>\varphi\left(L_{t}\right)\right\}$ is in the tail sigma
field, so its probability is $0$ or $1$. From Proposition
\ref{estimloc}, if
$\int_{0}^{\infty}\frac{dx}{\varphi^{2\mu}\left(x\right)}<\infty$,
then $\mathbb{P}\left(\exists t\geq
0,\;R_{t}>\varphi\left(L_{t}\right)\right)<1$, and hence
$\mathbb{P}\left(\forall A>0,\;\exists t\geq
A,\;R_{t}>\varphi\left(L_{t}\right)\right)=0$. Next, if
$\int_{0}^{\infty}\frac{dx}{\varphi^{2\mu}\left(x\right)}=\infty$,
then, again from Proposition \ref{estimloc}, we have, for all
$u\geq0$:
$$\mathbb{P}\left(\exists t\geq
\tau_{u},\;R_{t}>\varphi\left(L_{t}\right)\right)=1.$$To conclude,
it suffices to notice that $\left\{\forall A>0,\; \exists t\geq
A,\;R_{t}>\varphi\left(L_{t}\right)\right\}$ is the decreasing
limit, as $u\rightarrow\infty$, of the events $\left\{\exists t\geq
\tau_{u},\;R_{t}>\varphi\left(L_{t}\right)\right\}$.
\end{proof}
\begin{rem}
The results of Proposition \ref{estimloc} and its corollary can also
be obtained with the help of  excursion theory for the standard
Brownian Motion and then for the Bessel processes by  time change
techniques.
\end{rem}

\section*{Acknowledgements}
I would like to thank Marc Yor for many fruitful discussions and for
sharing with me some of his unpublished notes (\cite{donatimartin}).

\end{document}